\magnification 1200
\input amssym.def
\input amssym.tex
\parindent = 40 pt
\parskip = 12 pt
\font \heading = cmbx10 at 12 true pt
 at 22 true pt
\font \medheading =cmbx7 at 14 true pt
 at 7 true pt
\def \R{{\bf R}}
\def \C{{\bf C}}
\def \Q{{\bf Q}}

\centerline{\medheading Applications of an elementary resolution of singularities algorithm}
\centerline{\medheading to exponential sums and congruences modulo $p^n$}
\rm
\line{}
\line{}
\centerline{\heading Michael Greenblatt}
\line{}
\centerline{December 9, 2014}
\baselineskip = 12 pt
\font \heading = cmbx10 at 14 true pt
\line{}
\line{}
\noindent{\bf 1. Introduction and statement of results}

\vfootnote{}{This research was supported in part by NSF grants DMS-0919713 and DMS-1001070}\noindent 

In this paper, we use the resolution of singularities algorithm of [G4] to generalize  to arbitrary local fields of characteristic zero the theorems of [G3] on $\R^n$ sublevel set volumes and oscillatory integrals with real phase function. The proofs of these 
generalizations use various aspects of the resolution of singularities algorithms of [G4] (but for the most part not the actual resolution of singularities theorems themselves.) The $p$-adic cases of our results provide new estimates for exponential 
sums as well as new bounds on how often a function $f(x)$ such as a polynomial with integer 
coefficients is divisible by various powers of a prime $p$ when $x$ is an integer. Thus we use classical analysis resolution of 
singularities methods on a class of problems traditionally approached using toric resolution of singularities techniques.

The estimates of this paper for the sublevel set measures and oscillatory integral decay rates will as in [G3] be expressed 
in terms of properties of the Newton polyhedron of the phase function $f(x)$. Such estimates go back to [V] and for $p$-adic fields
there is an  extensive body of research on such estimates. We mention [D][DHo] [DLo][LM][Ve][Wr][Zu1][Zu2] for a sampling. Unlike
in many such papers, especially in higher dimensions, we will not require a nondegeneracy condition such as that of [V]. Instead, as in [G3]
 our theorems will be stated in terms of the orders of the zeroes of certain
polynomials $f_F(x)$ associated to $N(f)$, given in Definition 1.5 below, and will go beyond such a nondegeneracy condition. 
We will have results for any Newton distance (see Definiton 1.6),  and the 
Newton distance will determine the conditions required on the maximal order of the zeroes of the above polynomials. Generalizations  when the Newton distance is less than or equal to 1 will be discussed after the statement of Theorem 1.2.

\noindent In this paper, a real oscillatory integral denotes an expression
$$I(\lambda) = \int_{\R^n}e^{i\lambda f(x)}\phi(x)\,dx \eqno (1.1)$$
Here $\phi(x)$ is a cutoff function defined on an appropriately small bounded open set containing the origin. Using
resolution of singularities it can be shown that if $\nabla f(0) = 0$ (the nontrivial case),
 then  $I(\lambda)$ has an asymptotic
 expression as $\lambda \rightarrow \infty$ of the form
$$I(\lambda) = c_{\phi} e^{i\lambda f(0)}(\ln(\lambda))^m\lambda^{-\delta} + o((\ln(\lambda))^m\lambda^{-\delta}) \eqno (1.2)$$ 
Here $m$ is a nonnegative integer and $\delta >0$, both  independent of $\phi$, and given any sufficiently small neighborhood $U$ of the origin $c_{\phi}$ is nonzero for at least one $\phi$ supported in $U$.

To understand what the complex analogue of $(1.1)$ might be, note that a key characteristic of oscillatory integrals
$(1.1)$ is that for any $\lambda$ the function $e^{i\lambda t}$ is a (continuous) additive character on $\R$. Furthermore, all
continuous additive characters of $\R$ are of this form. Thus with the view of finding analogues to $(1.1)$ we are interested in
what the continuous additive characters of $\C$ are. Suppose $\chi(z)$ is one such character. Then $\chi(z) = \chi(Re(z),0)\chi(0,Im(z))
= e^{i\lambda_1 Re(z)}e^{i\lambda_2 Im(z)}$ for some $\lambda_1$ and $\lambda_2$. There necessarily exists some 
complex number $w$ such that for any $z$, $Re(wz) = \lambda_1Re(z) + \lambda_2Im(z)$. Thus the continuous additive characters
of $\C$ are the functions of the form $e^{iRe(wz)}$ for some $w \in \C$. Hence for our purposes natural analogues of $(1.1)$ 
are the integrals
$$I(w) = \int_{\C^n}e^{i Re(w f(z))}\phi(z)\,dz \eqno (1.3)$$
Here again $\phi(z)$ is a cutoff function, and the goal now will be to find optimal decay estimates 
$|I(w)| \leq C(\ln|w|)^m |w|^{-\delta}$ as $|w| \rightarrow \infty$. As in the real case, if $f(z)$ has some nonvanishing first derivative 
at the origin then one gets arbitrarily fast decay, so we always assume this is not the case here.

We now direct our attention to characters on $p$-adic fields. Any $p$-adic number $x$ may be written as 
$x = \sum_{k = k_0}^{\infty} b_k p^k$, where each $b_k \in \{0,...,p-1\}$, $b_{k_0} \neq 0$, and $|x| = p^{-k_0}$.
Addition on the $p$-adics in such a form is done as one adds natural numbers written in base $p$, using carrying. It is 
well-known (and relatively easy to show) that the 
continuous additive characters on the $p$-adics are functions of the form $\chi(x) = \xi(yx)$, where $y \in {\Q}_p$ and where
 $\xi( \sum_{k \geq k_0} b_k p^k) = e^{2\pi i ( \sum_{k = k_0}^{-1} b_k p^k)}$ for $k_0 < 0$ and $\xi(x) = 1$ for $k_0 \geq 0$. 

One can determine the additive characters on any $p$-adic field (a finite extension field of some $\Q_p$) from the characters
on $\Q_p$ analogously to how the additive characters on $\C$ were determined above from those on $\R$. Namely, suppose 
$K$ is a field extension of $\Q_p$ of degree $l$. Then we may write elements $x$ of $K$ in the form $(x_1,...,x_l)$ where each
$x_k \in \Q_p$. Thus if $\chi$ is a continuous additive character on $K$, one has
$$\chi(x) = \chi(x_1,0,...,0) \,...\,\chi(0,...,0,x_l)$$
By the form of the characters on $\Q_p$,  for some $y_1,...,y_l \in\ Q_p$ the above can be written as 
$$\chi(x) = \xi(y_1x_1)\,...\,\xi(y_lx_l) = \xi(y_1x_1 + ... + y_lx_l)$$
Analogously to the complex case, there is some $z \in K$ such that the first component of $zx$ is $y_1x_1 + ... + y_lx_l $ for
 all $x \in K$. Denoting this first component by $R(zx)$, for any $x$ one therefore has
$$\chi(x) = \xi(R(zx)) \eqno (1.4)$$
Thus a natural analogue of the oscillatory integral $(1.1)$ for $p$-adic fields is given by
$$I(z) = \int_{|x| < \delta}\xi(R(zf(x)))\,dx \eqno (1.5)$$
Since all smooth functions on a $p$-adic field are locally constant, instead of having a cutoff function $\phi(x)$ in $(1.5)$ we
restrict the domain of integration to ${|x| < \delta}$ for some $\delta$.

In the real case it is well known (see [AGV]) that sharp estimates for oscillatory integrals usually follow from sharp estimates for
the measure of sublevel sets; given a sufficiently small open set $U$ containing the origin one may look for the best possible estimate of the
form $|\{x \in U: |f(x)| < \epsilon\}| < C|\ln(\epsilon)|^m \epsilon^{\delta}$. This $(m,\delta)$ can then be translated into 
decay estimates for the oscillatory integral $(1.1)$, and this rate of decay is sharp except in certain
exceptional situations. This translation is proven using 
resolution of singularities to show both the sublevel set measures and the oscillatory integral decay have asymptotic expansions,
and then using integration by parts in a certain way to go from the sublevel set measures to the oscillatory integral 
estimates. The analgoue for $p$-adic fields was proved by Igusa [I1]-[I3], and in Theorem 2.2 we will prove the corresponding
statement for $K = \C$ using results from [G4].

\noindent We now give some relevant definitions.

\noindent {\bf Definition 1.4.} Let $f(x)$ be a function such that $f(x)$ has a convergent power series expansion
$\sum_{\alpha} f_{\alpha}x^{\alpha}$ on a neighborhood of the origin in $K^n$.
For any $\alpha$ for which $f_{\alpha} \neq 0$, let $Q_{\alpha}$ be the octant $\{t \in \R^n: 
t_i \geq \alpha_i$ for all $i\}$. Then the {\it Newton polyhedron} $N(f)$ of $f(x)$ is defined to be 
the convex hull of all $Q_{\alpha}$.  

A Newton polyhedron can contain faces of various dimensions in various configurations. 
These faces can be either compact or unbounded. In this paper, as in earlier work such as [G3] and [V], an 
important role is played by the following functions, defined for compact faces of the Newton polyhedron. A vertex
 is always considered to be a compact face of dimension zero.

\noindent {\bf Definition 1.5.} Suppose $F$ is a compact face of $N(f)$. Then
if $f(x) = \sum_{\alpha} f_{\alpha}x^{\alpha}$ denotes the Taylor expansion of $f$ like above, 
define $f_F(x) = \sum_{\alpha \in F} f_{\alpha}x^{\alpha}$.

\noindent The statements of several of our theorems will use the following terminology.

\noindent {\bf Definition 1.6.} Assume $N(f)$ is nonempty. Then the 
{\it Newton distance} $d(f)$ of $f(x)$ is defined to be $\inf \{t: (t,t,...,t,t) \in N(f)\}$.

\noindent {\bf Definition 1.7.} The {\it central face} of $N(f)$ is the face of $N(f)$ of minimal dimension intersecting the line
$t_1 = t_2 = ... = t_n$. 

\noindent In Definition 1.7, the central face of $N(f)$ is well-defined since it is given by the intersection of all faces of $N(f)$
intersecting the line $t_1 = t_2 = ... = t_n$. An equivalent definition that is sometimes used (such as in [AGV]) is that the central
face of $N(f)$ is the unique face of $N(f)$ intersecting the line $t_1 = t_2 = ... = t_n$ in its interior.

We now come to our theorems concerning sublevel set measures and oscillatory integrals. They are analogues of corresponding
results in [G3]. $f(x)$ will always denote a function with convergent power series on a neighborhood of the origin with $f(0) = 0$. The statements of the theorems are slightly different for different fields $K$, in that they will depend on the
 dimension of $K$ over its base field. Correspondingly, in the following $b_K = 1$
if $K = \R$, $b_K = 2$ if $K = \C$, and for an extension of the $p$-adics  $b_K$ denotes the degree of $K$ over $\Q_p$.

We use the notation $|A|$ to denote the measure of a set $A$. If $K = \R$ or $\C$ this denotes the usual Lebesgue measure.
When $K$ is ${\bf Q}_p$ we use the traditional Haar measure that assigns measure 1 to $\{x: |x| \leq 1\}$, and if $K$ is a
finite extension of ${\bf Q}_p$ we use the product measure induced by that of ${\bf Q}_p$. As for which valuation we use
on $K$, we will use the traditional $|x|_p = p^{-v_p(x)}$ valuation on ${\bf Q} \subset {\bf Q}_p$ and its natural extension
 to $K$ if $K$ is a 
finite extension of ${\bf Q}_p$. Note that the $b_K$th power of this valuation is used in [Cl] which accounts for the difference
 in the exponents in the theorems such as in section 2.

\noindent {\bf Theorem 1.1.} Let $K_0$ denote $K - \{0\}$. For a compact face $F$ of $N(f)$, let $o(F)$ denote the maximum 
order of any zero of $f_F(x)$ in $K_0^n$. Let $h$ denote the dimension of the central face of $N(f)$. For a small open set 
$U$ containing the origin, let $g(\epsilon)$ denote the measure of $\{x \in U: |f(x)| < \epsilon\}$. Then if 
$U$ is sufficiently small, there are positive constants $C$ and $C'$ depending on $U$ and $f$ such that the following hold
for $0 < \epsilon < {1 \over 2}$.

\noindent {\bf a)} If $o(F) \leq d(f)$ for all compact faces $F$ of $N(f)$, with $o(F) < d(f)$ when $F$ is a subset of the
central face of $N(f)$, then 
$$C|\ln(\epsilon)|^{n-h-1}\epsilon^{b_K \over d(f)} \leq g(\epsilon) \leq C'|\ln(\epsilon)|^{n-h-1}\epsilon^{b_K \over d(f)}$$
\noindent {\bf b)} If $o(F) \leq d(f)$ for all compact faces $F$ of $N(f)$ with $o(F) = d(f)$ for at least one compact $F$  
contained in the central face of $N(f)$, then
$$C|\ln(\epsilon)|^{n-h-1}\epsilon^{b_K \over d(f)} \leq g(\epsilon) \leq C'|\ln(\epsilon)|^{n-h}\epsilon^{b_K \over d(f)}$$
\noindent {\bf c)} If $o(F) > d(f)$ for at least one compact face $F$ of $N(f)$, let $s(f)$ denote $\sup_F o(F)$. Then we have
$$C|\ln(\epsilon)|^{n-h-1}\epsilon^{b_K \over d(f)} \leq g(\epsilon) \leq C'\epsilon^{b_K \over s(f)}$$

\noindent {\bf Theorem 1.2.} If in Theorem 1.1 one has an upper bound $g(\epsilon) \leq C'|\ln(\epsilon)|^m\epsilon^{\delta}$, 
then for sufficiently large $|\lambda|$, $|w|$, or $|z|$
the oscillatory integral  $(1.1)$, $(1.3)$, or $(1.5)$ respectively satisfies the analogous bound $|I(\lambda)|  \leq C
(\ln|\lambda|)^m |\lambda|^{-\delta}$,  $|I(w)| \leq C(\ln|w|)^m |w|^{-\delta}$, or 
$|I(z)| \leq C(\ln|z|)^m |z|^{-\delta}$.

In the case where $K = \R$, Theorems 1.1 and 1.2 were proven in [G3], where the oscillatory integral upper bounds corresponding
to parts a) and b) of Theorem 1.1 were shown typically to be
sharp. It is unclear what the sharpness situation is for oscillatory integrals in the $p$-adic case. The upper bounds of part c)
of Theorem 1.1 and the corresponding upper bounds of Theorem 1.2 are usually not sharp; to get sharp estimates one normally
needs more detailed information about the singularities of $f(x)$ than the Newton polyhedron provides. 

It is worth mentioning that 
there is an additional situation where Theorem 1.2  is known to hold, namely where $d(f) \leq 1$ 
and $o(F) \leq 1$ for all compact faces $F$ of $N(f)$. This was proved using toric resolution of singularities in [V], and the 
method contained therein generalizes to the $K$ considered in this paper. Since in this paper Theorem 1.2 is effectively a 
consequence of Theorem 1.1, where the analogous statement is false, this case is not covered by Theorem 1.2. However, it 
is possible to prove this case directly using Lemma 3.1 through a direct integration
by parts of the exponential since the phase function will have nonvanishing gradient.

Often in papers using the toric methods of [AGV], theorems explicitly make the assumption that
$$\{x \in U: \nabla f(x) = 0\} \subset \{x \in U: f(x) = 0\} \eqno (1.6)$$
One can use resolution of singularities to see that $(1.6)$ always holds in a sufficiently small neighborhood of the origin  (when $f(0) = 0$). Namely, let $\Psi(x)$ be as in Hironaka's theorem such that $f \circ \Psi(x)$ is locally comparable to a monomial; that is, in a local
coordinate system $f \circ \Psi(x)$ is of the form $a(x)m(x)$ where $m(x)$ is a monomial and $a(x)$ is nonvanishing. Then by 
the chain rule, $(1.6)$ holds on a neighborhood of the origin if the analogous statement for $f \circ \Psi$ holds on a neighborhood
of $\Psi^{-1}(0)$. The latter statement can be seen to be true by a direct computation in coordinates for which $f \circ \Psi(x)$ is of the this form $a(x)m(x)$. 

Thus while we don't explicitly assume $(1.6)$, we do only prove Theorems 1.1 and 1.2  on a neighborhood of the origin that is
sufficiently small for our various lemmas to work. Because of the involved nature of the arguments leading to these theorems, it is hard to ascertain if this neighborhood is also sufficiently small that $(1.6)$ necessarily holds throughout.

It is natural to ask to what extent the arguments giving Theorem 1.1 and 1.2 extend to local fields of positive characteristic. The argument  in
[G4] that proves Lemma 3.5 of that paper (which is the same as Lemma 3.1 is this paper) makes use of the fact that the field has characteristic zero; equation
$(3.1b)$ does not immediately hold in the positive characteristic case by the same argument. In addition, we utilize the Weierstrass preparation
theorem for characteristic zero in the proof of Theorem 2.1, so any extension to the positive characteristic case would require an
analogue or a substitute. So any adaptation of our proofs to the positive characteristic case would have to take into account both issues.

\noindent {\bf Number theoretic consequences.}

The cases $K = \Q_p$ of Theorems 1.1 and 1.2 have some number-theoretic consequences. In Theorems 1.3
and 1.4, $f(x)$ denotes a power series $\sum f_{\alpha} x^{\alpha}$ with integer coefficients 
that converges on a neighborhood of the origin when viewed as a power series on ${\bf Q}_p^n$, satisfying $f(0) = 0$.
The following theorem is a consequence of Theorem 1.1; the condition that each $a_i$ is sufficiently large ensures that
Theorem 1.1 applies to $f_a(x)$ on an entire ball of radius 1 centered at the origin.

\noindent {\bf Theorem 1.3.} For $a = (a_1,...,a_n)$, each $a_i$ an integer, let $f_a(x_1,...,x_n)$ denote the function 
$f(p^{a_1}x_1,...,p^{a_n}x_n)$. There exists a constant $M > 0$ depending on $f(x)$ such that if $a_i > M$ for all $i$ then
there are constants $C, C' > 0$ depending on $a, f(x),$ and $p$ such that the following hold for every positive integer $l$.

\noindent {\bf a)} If $o(F) \leq d(f)$ for all compact faces $F$ of $N(f)$, with $o(F) < d(f)$ when $F$ is a subset of the
central face of $N(f)$, then
$$ Cl^{n-h-1}p^{- {l \over d(f)}} \leq {1 \over p^{ln}}\,\,\,\#\{x \in \{0,...,p^l - 1\}^n: p^l {\rm\,\, divides\,\,} f_a(x) \} \leq C'l^{n-h-1}p^{- {l \over d(f)}}$$
\noindent {\bf b)}  If $o(F) \leq d(f)$ for all compact faces $F$ of $N(f)$ with $o(F) = d(f)$ for at least one compact $F$  
contained in the central face of $N(f)$, then 
$$ Cl^{n-h-1}p^{- {l \over d(f)}} \leq {1 \over p^{ln}}\,\,\,\#\{x \in \{0,...,p^l - 1\}^n: p^l {\rm\,\, divides\,\,} f_a(x)\} \leq C'l^{n-h}p^{- {l \over d(f)}}$$
\noindent {\bf c)} If $o(F) > d(f)$ for at least one compact face $F$ of $N(f)$, let $s(f)$ denote $\sup_F o(F)$. Then
$$ Cl^{n-h-1}p^{- {l \over d(f)}}  \leq {1 \over p^{ln}}\,\,\,\#\{x \in \{0,...,p^l - 1\}^n: p^l {\rm\,\, divides\,\,} f_a(x)\} \leq C'p^{- {l \over s(f)}}$$
Similarly, the oscillatory integral result Theorem 1.2 implies

\noindent {\bf Theorem 1.4.} Let $f_a(x_1,...,x_n)$ be as in Theorem 1.3. There exists a constant $M > 0$ depending on $f(x)$ such that if $a_i > M$ for all $i$ then
there is a constant $C > 0$ depending on $a, f(x),$ and $p$ such that the following hold for every positive integer $l$.

\noindent {\bf a)} If $o(F) \leq d(f)$ for all compact faces $F$ of $N(f)$, with $o(F) < d(f)$ when $F$ is a subset of the
central face of $N(f)$, then
$${1 \over p^{ln}}\,\,\,\bigg|\sum_{x \in   \{0,...,p^l - 1\}^n} e^{2\pi i{ f_a(x) \over p^l}}\bigg| \leq C l^{n-h-1}p^{ - {l \over d(f)}}$$
\noindent {\bf b)}  If $o(F) \leq d(f)$ for all compact faces $F$ of $N(f)$ with $o(F) = d(f)$ for at least one compact $F$  
contained in the central face of $N(f)$, then 
$${1 \over p^{ln}}\,\,\,\bigg|\sum_{x \in   \{0,...,p^l - 1\}^n} e^{2\pi i{ f_a(x) \over p^l}}\bigg| \leq C l^{n-h}p^{ -{ l \over d(f)}}$$
\noindent {\bf c)} If $o(F) > d(f)$ for at least one compact face $F$ of $N(f)$, let $s(f)$ denote $\sup_F o(F)$. Then
$${1 \over p^{ln}}\,\,\,\bigg|\sum_{x \in   \{0,...,p^l - 1\}^n} e^{2\pi i{ f_a(x) \over p^l}}\bigg| \leq C p^{ - {l \over s(f)}}$$
Since in Theorems 1.3 and 1.4 the constants $C$ and $C'$ do depend on $p$, no uniformity 
in $p$ is proven here such as was conjectured in the homogeneous case by Igusa [I3]. Also note that if $f(x^{b_1},...,
x^{b_n})$ is homogeneous for some positive integers $b_i$ then one can automatically replace $f_a(x)$ with $f(x)$ in
Theorems 1.3 and 1.4 by choosing $a$ appropriately (although the constants will change).

\noindent {\bf 2. Van der Corput lemmas}

As in a number of papers that give asymptotics for oscillatory integrals and related matters, we will make significant use of
lemmas related to the classical Van der Corput lemma. When $K = \R$ what we will need will readily follow from the classical
one-dimensional Van der Corput lemma:

\noindent {\bf Classical Van der Corput Lemma.} Suppose $g(x)$ is a $k$ times differentiable function on an interval $I$ 
satisfying  $|g^{(k)}(t)| > \eta > 0$ for all $t$. Then for a constant $A_k$ depending only on $k$, for any $\epsilon > 0$ we have
$$|\{t \in I: |g(t)| < \epsilon\}| < A_k \epsilon^{1 \over k}\eta^{-{1 \over k}} \eqno (2.1a)$$
If $\phi(x)$ is a $C^1$ function on $I$ and $k > 1$, then for some constant $B_{\phi,k}$ one has
$$\bigg|\int_I e^{i\lambda g(t)} \phi(t)\,dt\bigg| < B_{\phi,k}(1 + |\lambda|)^{-{1 \over k}} \eqno (2.1b)$$
If $g'(t)$ is piecewise monotone on $I$, then $(2.1b)$ also holds for $k = 1$, where the constant will now
also depend on the number of pieces on which $g'(t)$ is monotone.

For $K$ other than $\R$, we will also make use of versions of the Van der Corput lemma that hold for analytic functions. The
sublevel set version we will use can be stated as follows.

\noindent {\bf Theorem 2.1.} Suppose $U \subset K^n$ is a bounded open set with $0 \in U$, and suppose $f(x)$ is a function whose Taylor
series at the origin converges on a neighborhood of $cl(U)$. Suppose also that there is some $k > 0$ such that
$|\partial_{x_n}^k f(x)| \neq 0$ on $cl(U)$. Then there is a constant $B$ independent of $\epsilon$ (but depending on $f(x)$)
such that for all $\epsilon > 0$ one has
$$|\{x \in U: |f(x)| < \epsilon\}| < B\epsilon^{b_K \over k} \eqno (2.2)$$

\noindent {\bf Proof.} The case where $K = \R$ follows immediately by localizing and then applying the Van der Corput Lemma 
in the $x_n$ direction. Suppose now $K = \C$. It suffices to prove $(2.2)$ in a neighborhood of any $x_0 \in cl(U)$. If $f(x_0) 
\neq 0$, this is immediate. So we suppose $f(x_0) = 0$. Let $l > 0$ be minimal such that $\partial_{x_n}^l f(x_0) = 0$; by assumption, $l \leq k$. Let $F(x) = f(x_0 + x)$. By the Weierstrass preparation theorem, there is an open set
$V$ containing $0$ on which we may write
$$F(x) = c(x)\big(x_n^l + \sum_{i = 0}^{l-1} a_i(x_1,...,x_{n-1})x_n^i\big) \eqno (2.3)$$
Here $c(x) $ is analytic with $c(0) \neq 0$. Since $c(0) \neq 0$, shrinking $V$ if necessary it suffices to show that for a 
constant $C$ independent of $\epsilon$, for all $\epsilon > 0$ we have
$$\big|\{x \in V: |x_n^l + \sum_{i = 0}^{l-1} a_i(x_1,...,x_{n-1})x_n^i| < \epsilon \}\big|< C\epsilon^{ 2 \over k} \eqno (2.4)$$
By the fundamental theorem of algebra we may factorize
$$x_n^l + \sum_{i = 0}^{l-1} a_i(x_1,...,x_{n-1})x_n^i = \prod_{i=1}^{l}(x_n - \alpha_i(x_1,...,x_{n-1}))\eqno (2.5)$$
In order for the product on the right of $(2.5)$ to be less than $\epsilon$ in magnitude, at least one $|x_n - \alpha_i(x_1,...,x_{n-1})|$ must be less than $\epsilon^{1 \over l}$. Hence given $(x_1,...,x_{n-1})$, the set of 
$x_n$ for which the magnitude of the product is less than $\epsilon$ has measure at most $l \epsilon^{2 \over l}
\leq k \epsilon^{2 \over k}$. Integrating this over all $(x_1,...,x_{n-1}) \in V$ gives $(2.4)$ as needed.
This completes the proof for $K = \C$.

If $K$ is a $p$-adic field, we do a similar argument using the $p$-adic Weierstrass preparation theorem (see Theorem 6.2.10 of
[Go]). As in the complex case, we may focus our attention on an $x_0 \in cl(U)$ for which $f(x_0) = 0$, and we define
$F(x) = F(x_0 + x)$. Again let $l$ be minimal such that $\partial_{x_n}^l f(x_0) = 0$. By the $p$-adic Weierstrass 
preparation theorem, there are balls $B_1 \subset K$ and $B_2
\subset K^{n-1}$ centered at the origins of $K$ and $K^{n-1}$ respectively such that if $(x_1,...,x_{n-1}) \in B_2$, 
then there are $b_0(x_1,...,x_{n-1}),...,b_{l-1}(x_1,...,x_{n-1}) \in K$ and a constant $A > 0$ such that for $x_n \in B_1$ one has
$$|F(x)| = A|x_n^l + \sum_{i=0}^{l-1}b_i(x_1,...,x_{n-1})x_n^i|\eqno (2.6)$$
Hence given $(x_1,...,x_{n-1}) \in B_2$, there is a finite extension $L$ of $K$ and $\beta_i(x_1,...,x_{n-1}) \in L$ such that
$$|F(x)|  = A|\prod_{i = 1}^l(x_n - \beta_i(x_1,...,x_{n-1}))| \eqno (2.7)$$
Thus if $|F(x)| < \epsilon$, there is some $i$ for which $|x_n - \beta_i(x_1,...,x_{n-1})| < A^{-{1 \over l}}\epsilon^{1 \over l}$ in $L$.
By the ultrametric property of $p$-adic fields, if there is any $y \in K$ such that $|y - \beta_i(x_1,...,x_{n-1})|  <A^{-{1 \over l}}
\epsilon^{1 \over l}$, then for all $z \in K$ with  $|z - \beta_i(x_1,...,x_{n-1})|  < A^{-{1 \over l}}\epsilon^{1 \over l}$ one has 
$|z - y| < A^{-{1 \over l}}\epsilon^{1 \over l}$ as well. Hence either the set of $x_n \in B_1$ for which $|x_n - \beta_i(x_1,...,x_{n-1})| < A^{-{1 \over l}}\epsilon^{1 \over l}$ is empty, or it is a subset of the set of $x_n \in B_1$ for which $|x_n - y_i(x_1,...,x_{n-1})| < A^{-{1 \over l}}\epsilon^{1 \over l}$, 
where $y_i(x_1,...,x_{n-1})$ is some element of $K$. This set has measure bounded by $A^{-b_K \over l}\epsilon^{b_K \over l}$. Hence for our fixed $(x_1,...,x_{n-1})$, the set of $x_n \in B_1$ with $|F(x)| < \epsilon$ has measure at most $ A^{-b_K \over l}l
 \epsilon^{b_K \over l} \leq A^{-b_K \over l} k \epsilon^{b_K \over k}$. Like with $K = \C$, we now integrate this over all $(x_1,...,x_{n-1}) \in B_2$, obtaining
$(2.2)$ as needed. This completes the proof of Theorem 2.1.

For the oscillatory integrals $(1.1)$, $(1.3)$, $(1.5)$ the analogue to Theorem 2.1 holds as well. For $K = \R$ this is once again an 
immediate consequence of the Van der Corput lemma. When $K$ is a $p$-adic field, it is a consequence of a recent result of
Cluckers [Cl]; see also [R] for an earlier partial result. The proof in [Cl] appears to extend also to $K = \C$. For the purposes of 
our paper however we need the following result which relies on resolution of singularities:

\noindent {\bf Theorem 2.2}. Suppose $f(x) = \sum_{\alpha} f_{\alpha}x^{\alpha}$ converges on a neighborhood of the origin in
$K^n$ with $f(0) = 0$. Suppose $\delta$ and $l$ are such that for any sufficiently small open $U$ containing the origin, 
for all $0 < \epsilon < {1 \over 2}$ one has 
$$|\{x \in U: |f(x)| < \epsilon\}| < A_{f,U}\,\epsilon^{\delta}|\ln(\epsilon)|^l\eqno (2.8)$$
Then there is a neighborhood $V$ of the origin such that if the support of the integrand of the oscillatory integral $(1.1)$,
$(1.3)$, or $(1.5)$ is contained in $V$, then for sufficiently large $|\lambda|$, $|w|$, or $|z|$ respectively
the oscillatory integral satisfies the analogous bound $|I(\lambda)|  \leq C
(\ln|\lambda|)^l |\lambda|^{-\delta}$,  $|I(w)| \leq C(\ln|w|)^l |w|^{-\delta}$, or 
$|I(z)| \leq C(\ln|z|)^l |z|^{-\delta}$.

\noindent {\bf Proof.} When $K = \R$, Theorem 2.2 is a relatively straightforward consequence of the existence of asymptotic 
expansions for sublevel set measures and oscillatory integrals; we refer to [AGV] Ch 7 for details. An elementary proof
based on an earlier version of the resolution of singularities algorithm of [G4] was also given in [G2]. When $K$ is a $p$-adic 
field, Theorem 2.2 was shown by Igusa [I1]-[I3]. 

It remains to consider $K = \C$. We will prove Theorem 2.2 in the complex case using results from [G4] (and it can be proved in a very similar fashion using Hironaka's theorem [Hi1]-[Hi2]); the proof will
also be quite similar to aforementioned argument of [I1]-[I3] for the $p$-adic field
case. So suppose $K = \C$ now. By Theorem 1.1 of [G4], there are a finite collection of functions $\beta_m(z)$
such that if $\phi(x)$ is a cutoff function supported on a sufficiently small neighborhood of the origin, one can write $\phi(x) = \sum_{m=1}^N \tilde{\phi}_m(x)$
such that each $\phi_m(z) = \tilde{\phi} \circ \beta_m(z)$ can be changed on a set of measure zero to become a smooth function. Furthermore
each $\beta_m(z)$ is one-to-one almost everywhere and on the support of $\phi_m$, and both 
$f \circ \beta_m(z)$ and the Jacobian $J_m(z)$ of $\beta_m(z)$ are of the form $a_m(z)p_m(z)$ with $a_m(z)$ nonvanishing and $p_m(z)$ a monomial. So for any measurable bounded function $h(x)$ we have
$$\int h(x)\phi(x)\,dx = \sum_{m=1}^N  \int h \circ \beta_m(z) |J_m(z)|^2\phi_m(z)\,dz \eqno (2.9)$$
Here $\phi_m(z)$ is smooth and on the support of the integrand of a given term of
the right-hand side of $(2.9)$ the function 
$f \circ \beta_m(z)$ is of the form $a_m(z)p_m(z)$  where $a_m(z)$ is nonvanishing and $p_m(z)$ is a nonconstant 
monomial. Applying $(2.9)$ to $h(x) = e^{i Re(w f(x))}$ for $w \in \C$ we get
$$\int  e^{i Re(w f(x))} \phi(x)\,dx = \sum_{m=1}^N \int e^{i Re(w a_m(z)p_m(z))} |J_m(z)|^2\phi_m(z)\,dz \eqno (2.10)$$
Since we may restrict ourselves to $x$ in an arbitrarily 
small neighborhood of the origin, in a given term of $(2.10)$ we may restrict ourselves to $z$ in an arbitrarily small neighborhood of $\beta_m^{-1}(0)$. Write $p_m(z) = \prod_{i = 1}^n z_i^{q_i}$. If $z \in \beta_m^{-1}(0)$, then since
$f(0) = 0$ one has $f \circ \beta_m (z) = f(0) = 0$ and therefore $p_m(z) = 0$. As a result, $\beta_m^{-1}(0)$ is a 
subset of $\{z \in \C^n: z_i = 0$ for some $i$ with $q_i > 0\}$. Using a partition of unity, splitting $\phi_m$ into a finite sum of bump functions with smaller support as needed, without loss of generality we can assume for some $z' \in \beta_m^{-1}(0)$ that $\phi_m$ is supported on a small neighborhood of $z'$ such that for some $i'$  with $z_i' = 0$ and $q_i > 0$, we can do a change of variables 
 such that what was once $a_m(z)z_i^{q_i}$ becomes $z_i^{q_i}$. (Recall that
$a_m(z)$ is nonvanishing). Thus in place of a given term of $(2.10)$ we may focus on
$$ \int e^{i Re(w\prod_{i=1}^n z_i^{q_i})} |J_m(z)|^2\phi_m(z)\,dz \eqno (2.11)$$
Note that $\phi_m(z)$ may be a different bump function from that of $(2.10)$ due to the localization and variable change.
Changing variable names if necessary we can assume $q_1 > 0$ in $(2.11)$. The idea now will be to use real
integrations by parts in the $z_1$ variable in $(2.11)$ and then integrate the result in the remaining variables. To
this end, for fixed $(z_2,...,z_n)$ we do a variable change of the form $z_1 \rightarrow e^{i\theta}z_1$ so that the integral 
of the integrand of
$(2.11)$ in the $z_1$ variable becomes of the form
$$ \int e^{i Re(|w|(\prod_{i=2}^n |z_i|^{q_i})z_1^{q_1})} |J_m(z)|^2\tilde{\phi}_m(z)\,dz_1\eqno (2.12a)$$
We divide this dyadically as
$$\sum_j \int e^{i Re(|w|(\prod_{i=2}^n |z_i|^{q_i})z_1^{q_1})} |J_m(z)|^2\phi_{mj}(z)\,dz_1\eqno (2.12b)$$
Here $\phi_{mj}(z)$ is supported in $2^{-j-1} < |z_1| < 2^{-j + 1}$.
Observe that for fixed $(z_2,...,z_n)$ the phase in $(2.12b)$ is of the form $A Re(z_1^{q_1})$, which is a homogeneous
polynomial of degree $q_1$ with  gradient of magnitude $\sim A|z_1|^{q_1 - 1}$. Hence one can integrate by parts as many times as one
wants; each time one gains a factor of ${C \over A|z_1|^{q_1- 1}}$ but one also loses a factor of ${C' \over |z_1|}$ each 
time a derivative  lands on $|J_m(z)|^2$ or $\phi_{mj}(z)$. Hence the net effect of $N$ integrations by parts is a factor
 of ${C_N \over A^N|z_1|^{q_1N}}$. As a result, the $j$th term of $(2.12b)$ is bounded by 
$$C_N\int_{(z_1,...,z_n) \in supp(\phi_{mj})} (|w|\prod_{i=1}^n |z_i|^{q_i})^{-N}|J_m(z)|^2\,d {z_1} \eqno (2.13)$$
By taking absolute values of the integrand and then integrating, the $j$th term of $(2.12b)$ is also bounded by
$$C \int_{(z_1,...,z_n) \in supp(\phi_{mj})}|J_m(z)|^2\,d {z_1} \eqno (2.14)$$
Combining $(2.13)$ and $(2.14)$, and then adding the result over all $j$, we see that $(2.12a)$ is bounded by
$$C_N'\int_{(z_1,...,z_n) \in supp(\phi_m)}\min (1,   (|w|\prod_{i=1}^n |z_i|^{q_i})^{-N})|J_m(z)|^2\,dz_1 \eqno (2.15)$$
Integrating this bound in the $z_2,...,z_n$ variables, we get that $(2.11)$ is bounded by
$$C_N'\int_{supp(\phi_m)}\min (1,   (|w|\prod_{i=1}^n |z_i|^{q_i})^{-N})|J_m(z)|^2\,dz_1\,...\,dz_n \eqno (2.16)$$
As long as $N$ is sufficiently large, by Lemma 3.2a) of [G3] for example, the integral $(2.16)$ will be bounded by a 
constant times the integral over the portion of the domain where $|w|\prod_{i=1}^n |z_i|^{q_i} < 1$, in other words the portion
where  $\prod_{i=1}^n |z_i|^{q_i} < {1 \over |w|}$. Since $f \circ \beta_m(z)$ is within a constant factor of 
$\prod_{i=1}^n |z_i|^{q_i}$ on the support of $\phi_m(z)$, if $N$ is taken large enough $(2.16)$ is at most
$$C_N'' \int_{\{z \in supp(\phi_m): \, f \circ \beta_m(z) < {1 \over |w|}\}}|J_m(z)|^2\,dz\eqno (2.17)$$
If we let $h(x)$ be the characteristic function of the set where $|f(x)| < {1 \over |w|}$ in  $(2.9)$, we see by adding  $(2.17)$ 
over all $m$ that $(2.11)$ is at most a constant times the measure of the portion of the support of $\phi(x)$ where $|f(x)| < {1 \over |w|}$. By $(2.8)$ this is at most $A_U(1 + |w|)^{-\delta}\ln(1 + |w|)^l$ and we are done with the
proof.

Note that although Theorems 2.1 and 2.2 taken together imply the analogue of Theorem 2.1 for oscillatory integrals,
one does not get a uniform constant this way as in the classical Van der Corput lemma. So in particular we do not recover the 
whole $p$-adic Van der Corput lemma of [Cl].

\noindent {\bf 3. Proofs of Theorems 1.1, 1.2, 1.3, and 1.4.}

\noindent We will make use of the following lemma from [G4].

\noindent {\bf Lemma 3.1.} (Lemma 3.5 of [G4]). Let $K_0 = K - \{0\}$ as before. For $0 \leq i \leq n-1$, let $\{F_{ij}\}$ denote the faces of $N(f)$  of dimension $i$. Let $o_{ij}$ denote the maximum order of the zeroes
of $f_{F_{ij}}(x)$ on $K_0^n$ (see Definition 1.5). Then for each $i$ and $j$ there are finitely many maps $\gamma_{ijk}(x)$, each a composition of blowups (and therefore a bijective monomial map when restricted to $K_0^n$), such that the following holds.

There is an open $E$ containing the origin and constants $a , b, b'> 0$
such that if $\phi(x) \in C_c(E)$, outside a set of measure zero one can write $\phi(x) = \sum_{ijk}\phi_{ijk}(x)$, where each $\phi_{ijk}$ is supported in $E$ and each $\rho_{ijk} = \phi_{ijk} \circ \gamma_{ijk}$ is a function on $K_0^n$ 
extending to a smooth function on all of $K^n$. This can be done in such a way that the following hold.

\noindent {\bf a)} Suppose $i > 0$. Let $B(0,a)$ denote  $\{y: |y_l| < a$ for all $1 \leq p\}$ and let
$Z_{ijk}= \{z: a < |z_m| < 2$ for all $1 \leq m \leq n - p\}$. Then for some $p$ there exists an open $Z_{ijk}' \subset Z_{ijk}$ such that 
$$  B(0,a) \times Z_{ijk}' \subset supp(\rho_{ijk}) \subset B(0,2) \times Z_{ijk}  \eqno (3.1a)$$
 Furthermore,  if $m(x)$ denotes $x^v$ for any vertex $v$ of $N(f)$ on $F_{ij}$, then there exists 
a directional derivative $\sum_m \beta_m \partial_{z_m}$ with $\sum_m |\beta_m| = 1$ and and an integer $0 \leq q \leq o_{ij}$ such that on $supp(\rho_{ijk})$ one has
$$b |m \circ \gamma_{ijk}(y,z)| \leq |(\sum_m \beta_m \partial_{z_m})^q (f\circ \gamma_{ijk})(y,z)| \leq b' |m \circ \gamma_{ijk}(y,z)| \eqno (3.1b)$$
\noindent {\bf b)} Suppose $i = 0$, so that $F_{0j}$ is a single vertex $v$. Then we have
$$B(0,a) \subset supp(\rho_{0jk}) \subset B(0,2) \eqno (3.2a)$$
If $m(x) = x^v$ then on $supp(\rho_{0jk})$ one has
$$ b |m \circ \gamma_{0jk}(y)| \leq | f \circ \gamma_{0jk}(y)| \leq b'|m \circ \gamma_{0jk}(y)|  \eqno (3.2b)$$
---------------------------------------------------------------------------------------------------------------------

Let $\phi(x)$ be a cutoff function satisfying the hypotheses of Lemma 3.1. Then by the change of variables formulas (which
holds for any $K$), if $j(x)$ is any bounded measurable function defined on the support of $\phi(x)$ we have
$$\int j(x)\phi(x)\,dx = \sum_{ijk} \int j(x)\phi_{ijk}(x)\,dx$$
$$= \sum_{ijk} \int j \circ \gamma_{ijk}(x) |J_{ijk}(x)|^{b_K}\rho_{ijk}(x)\,dx \eqno (3.3)$$
Here $J_{ijk}(x)$ denotes the Jacobian determinant of $\gamma_{ijk}(x)$, which is a monomial since each component of
$\gamma_{ijk}(x)$ is. 

By Lemma 3.3 of [G4], for a given $(i,j,k)$ with $i > 0$ there is a single $(a_1,...,a_p)$ such that if $m(x) = x^v$ for any vertex $v$ of $N(f)$ on $F_{ij}$, then $m \circ \gamma_{ijk}(x)$ is of the form $y_1^{a_1}...y_p^{a_p}z_1^{b_1}...z_{n-p}^{b_{n-p}}$ for some $(b_1,...,b_{n-p})$ which depends on $v$. When $i = 0$, we define the $a_l$'s by
 $y_1^{a_1}...y_n^{a_n} = m \circ \gamma_{0jk}(x)$, where $m(x) = x^v$ for $v$ the vertex of $N(f)$ corresponding to
$F_{0j}$. The next lemma
relates the $a_l$'s to $J_{ijk}(x)$ and the Newton distance $d(f)$.

\noindent {\bf Lemma 3.2.} Let $h$ be the dimension of the central face of $N(f)$ as in the statement of Theorems 1.1-1.4. 
Write the Jacobian determinant $J_{ijk}(x)$ of $\gamma_{ijk}(x)$ as $y_1^{e_1}...y_p^{e_p}z_1^{f_1}...z_{n-p}^{f_{n-p}}$ (with only $y$ variables if $i = 0$). Then for all $1 \leq l \leq p$ one has
$${a_l \over e_l+ 1} \leq d(f) \eqno (3.4)$$
Equality holds for at most $n - h$ values of $l$, and there is at least one $(i,j,k)$ such that equality does 
hold for $n - h$ values of $l$. Equality can only hold if $F_{ij}$ is a subset of the central face of $N(f)$.

\noindent {\bf Proof.} In [G4], the set of all $(a_1,...,a_{p})$ and $(e_1,...,e_{p})$ was determined by $N(f)$ and not the particular field $K$. Hence without loss of generality we may take $K = \R$. If one does
the coordinate change $y_l = Y_l^{1 \over e_l + 1}$ and $z_m = Z_m^{1 \over f_m + 1}$ for all $l,m$ and let $\Gamma_{ijk}(Y,Z)$ be the composition of
 $\gamma_{ijk}(x)$ with this coordinate change, then $\Gamma_{ijk}(Y,Z)$ has constant Jacobian determinant  and furthermore
if $m$ is as in Lemma 3.1 we have
$$m \circ \Gamma_{ijk}(Y,Z) = Y_1^{{a_1 \over e_1 + 1}}...Y_p^{{a_p \over e_p + 1}} Z_1^{b_1 \over f_1 + 1}...Z_{n-p}^
{b_{n-p} \over f_{n-p} + 1} \eqno (3.5)$$
The lemma is now an immediate consequence of Lemma 2.6 of [G3].

 In order to prove Theorem 1.1, we apply $(3.3)$ with $j(x)$ 
 the characteristic function
 of $\{x \in U: |f(x)| < \epsilon\}$ for a sufficiently small neighborhood $U$ of the origin. We will bound
 a given term $I_{ijk}$ of the sum $(3.3)$ for this $j(x)$ and add the bounds. This will give Theorem 1.1. We first consider a term $I_{ijk}$ such 
that $i = 0$ or such that $i > 0$ but $q = 0$ in
 $(3.1b)$. In these situations, since the $|z_m|$ are bounded above and below, by $(3.1b)-(3.2b)$ there are constants $C,C' > 0$ such that on the domain of integration of $I_{ijk}$ one has
$$C|y_1^{a_1}...y_p^{a_p}| \leq |f \circ \gamma_{ijk}(y,z)| \leq C'|y_1^{a_1}...y_p^{a_p}| \eqno (3.6a)$$
Similarly, the fact that the $|z_m|$ are bounded above and below imply that on the domain of integration of $I_{ijk}$ one has
$$C''|y_1^{e_1}...y_{p}^{e_p}| \leq |J_{ijk}(y,z)| \leq C'''|y_1^{e_1}...y_{p}^{e_p}| \eqno (3.6b)$$
In view of the inclusions $(3.1a)-(3.2a)$, $I_{ijk}$ therefore satisfies
$$C_1\int_{\{y \in B(0,a): |y_1^{a_1}..y_p^{a_p}| < {\epsilon \over C'}\} } |y_1^{e_1}...y_{p}^{e_p}|^{b_K} \,dy \leq I_{ijk}$$
$$I_{ijk} \leq 
C_2 \int_{\{y \in B(0,2): |y_1^{a_1}..y_p^{a_p}| < {\epsilon \over C}\}} |y_1^{e_1}...y_{p}^{e_p}|^{b_K}\,dy \eqno (3.7)$$
Suppose $K = \R$ now (so $b_K = 1$), and change variables in $(3.7)$, with the
new $y_j$ being the former $y_j^{e_j + 1}$. Then for some $a',b',c',$ and $d'$, $(3.7)$ becomes
$$C_3 \int_{\{y \in B(0,a'): |y_1^{a_1 \over e_1 + 1}...y_{p}^{a_{p} \over e_{p} + 1}| <  c'\epsilon\}} 1\,dy
 \leq I_{ijk} \leq C_4 \int_{\{y \in B(0,b'):|y_1^{a_1 \over e_1 + 1}...y_{p}^{a_p \over e_p + 1}|<  d'\epsilon\}} 1 \,dy \eqno (3.8)$$
Measures of sets of the form $\{y \in B(0,1): |y_1^{a_1 \over e_1 + 1}...y_{p}^{a_p \over e_p + 1}| < \delta\}$ 
were analyzed in [G3]. By Lemma 3.1a) of [G3], if $M_{ijk}$ denotes the measure of this set, $\rho$ denotes the maximum of the
${a_l \over e_l + 1}$, and $\sigma$ denotes  the number of times this value of ${a_l \over e_l + 1}$ occurs, we have
$$C_5\,\,\delta^{1 \over \rho}|\ln(\delta)|^{\sigma - 1} \leq M_{ijk} \leq C_6\,\,\delta^{1 \over \rho}|\ln(\delta)|^{\sigma - 1} \eqno (3.9)$$ 
Thus using $(3.9)$ and rescaling the balls of the left and right-hand sides of $(3.8)$, we get
$$C_7\,\,\epsilon^{1 \over \rho}|\ln(\epsilon)|^{\sigma - 1} \leq I_{ijk} \leq C_8\,\,\epsilon^{1 \over \rho}|\ln(\epsilon)|^{\sigma - 1} \eqno (3.10)$$ 
In the case where $K$ is a field other than $\R$, one can prove the analogous statement to $(3.9)$ as in the proof of Lemma 3.1a)
 of [G3], inducting on the dimension and using Fubini's theorem. The analogue one obtains is
$$C_5'\,\,\delta^{b_K\over \rho}|\ln(\delta)|^{\sigma - 1} \leq M_{ijk} \leq C_6'\,\,\delta^{b_K\over \rho}|\ln(\delta)|^{\sigma - 1} \eqno (3.9')$$
Thus we get the following generalization of $(3.10)$.
$$C_7'\,\,\epsilon^{b_K \over \rho}|\ln(\epsilon)|^{\sigma - 1} \leq I_{ijk} \leq C_8'\,\,\epsilon^{b_K \over \rho}|\ln(\epsilon)|^{\sigma - 1} \eqno (3.10')$$ 
Since by Lemma 3.2, ${1 \over \rho} \geq {1 \over d(f)}$ with $\sigma \leq n - h $ when $\rho = d(f)$, equation $(3.10')$ gives the required upper bounds for Theorem 1.1 for any $I_{ijk}$ with $i = 0$, or with $i > 0$ and $q = 0$.

We now consider the terms $I_{ijk}$ of the sum $(3.3)$ for $i > 0$ such that $q > 0$ in $(3.1b)$. Doing a linear change of
variables in the $z$ variables if necessary, without loss of generality we may assume that the directional derivative  $\sum_m \beta_m \partial_{z_m}$ is just $\partial_{z_n}$.
To find the needed upper bounds, we integrate the integrand of $I_{ijk}$ in the $z_n$ direction first, using 
Theorem 2.1, and then integrate the result in the remaining directions. Using $(3.1b)$ in conjunction with Theorem 2.1, we see
that this
one-dimensional integral in the $z_n$ direction is bounded by $C\epsilon^{b_K \over q}|y_1|^{-{a_1  b_K \over q}}...|y_p|^{-{a_p  b_K \over q}}|z_1|^{-{b_1b_K \over q}}...|z_{n-p-1}|^{-{b_{n-p-1}b_K \over q}}$. Since the $|z_m|$ are bounded away from 0, this is bounded by $C'\epsilon^{b_K \over q}|y_1|^{-{a_1  b_K \over q}}...|y_p|^{-{a_p  b_K \over q}}$. This one dimensional
integral is also bounded by a uniform constant since the domain of integration is 
bounded. Integrating the minimum of these two bounds  in the $(z_1,...,z_{n-1})$ variables, we obtain (using that $|J_{ijk}(y,z)| \sim |y_1^{e_1}...y_p^{e_p}|$ since the $|z_m|$ are bounded above and below)
$$I_{ijk} \leq C_1 \int_{B(0,2)} |y_1^{e_1}...y_p^{e_p}|^{b_K} \min(1, \epsilon^{b_K \over q}|y_1|^{-{a_1 b_K\over q}}...|y_p|
^{-{a_p b_K \over q}})\,dy \eqno (3.11)$$
Again we first focus on the $K = \R$ case. Changing variables from $y_j $ to $y_j^{e_j + 1}$, $(3.11)$ becomes
$$I_{ijk}  \leq C_2 \int_{B(0,2)}  \min(1, \epsilon^{1\over q} |y_1|^{-{a_1 \over q(e_1 + 1)}}...|y_p|
^{-{a_p \over q(e_p+ 1)}})\,dy \eqno (3.12)$$
Integrals of the form $(3.12)$ were also analyzed in [G3]. Again let $\rho = \max_j {a_j \over e_j + 1}$ and $\sigma$ 
the number of times this maximum is achieved. By Lemma 3.1 a) and d) of [G3], if $\rho > q$, one has
$$I_{ijk} \leq C_3 \epsilon^{1 \over \rho}|\ln(\epsilon)|^{\sigma - 1}\eqno (3.13a)$$
By Lemma 3.2 c) of [G3], if $\rho = q$ we have
$$I_{ijk} \leq C_4 \epsilon^{1 \over \rho}|\ln(\epsilon)|^{\sigma}\eqno (3.13b)$$
While by Lemma 3.2b) of [G3], if $\rho < q$ then we have
$$I_{ijk} \leq C_5 \epsilon^{1 \over q} \eqno (3.13c)$$
Note that the exponent of $\epsilon$ of $(3.13a)-(3.13c)$ can be succinctly written as $\min({1 \over q}, {1 \over \rho})$.
In the case where $K$ is not $\R$, one can prove inequalities analogous to $(3.13a)-(3.13c)$ similarly to above; one may bound
the right-hand side of $(3.11)$ to obtain analogues to $(3.13a)-(3.13c)$ by inducting on the 
dimension and using Fubini's theorem. The result is the same as $(3.13a)-(3.13c)$, except the exponents of $\epsilon$ are multiplied by $b_K$. In other words, the exponent becomes $\min({b_K  \over q}, {b_K \over \rho})$.
Note that this exponent is decreasing with increasing $q$. Thus if $s(f)$ denotes the maximum value of $q$ over all
faces of $N(f)$,
this exponent is at least $\min({b_K \over s(f)}, {b_K \over \rho})$. Furthermore by Lemma 3.2, 
$\rho \leq d(f)$, so that the exponent is at least $\min({b_K \over s(f)}, {b_K \over d(f)})$. Also, by Lemma 3.2 if 
$\rho = d(f)$ then $\sigma \leq n - h$, where $h$ is the dimension of the central face of $N(f)$. Equality may hold only
when $F_{ij}$ is a subset of the central face of $N(f)$. Hence when $s(f) < d(f)$ we have
$$I_{ijk} \leq C_6\epsilon^{b_K \over d(f)}|\ln(\epsilon)|^{n-h-1}\eqno (3.14a)$$
When $s(f) = d(f)$ we get
$$I_{ijk} \leq C_7 \epsilon^{b_K \over d(f)}|\ln(\epsilon)|^{n-h}\eqno (3.14b)$$
Here equality can only hold if $F_{ij}$ is a subset of the central face. When $s(f) > d(f)$ we get
$$I_{ijk} \leq C_8 \epsilon^{b_K \over s(f)} \eqno (3.14c)$$
These are exactly the exponents of Theorem 1.1. So adding over all $i,j,$  and $k$ proves all upper bounds of
 Theorem 1.1. Now Theorem 1.2 follows immediately via Theorem 2.2.

The lower bounds of Theorem 1.1 hold for the following reason. By Lemma 3.2, there is at least 
one face $F_{ij}$ for which $\rho = d(f)$ and $\sigma = n-h$.  The constructions of [G4] are such that for $F_{ij}$ one can ensure
that there is at least one $k$ for which $(3.2b)$ holds or for which $(3.1b)$ holds with $q= 0$.
As a result, the lower bounds of $(3.10')$ imply the required lower bounds of Theorem 1.1, 
completing the proof of that theorem.

\noindent {\bf Proof of Theorems 1.3 and 1.4.}

Theorems 1.3 and 1.4 are translations of Theorems 1.1 and 1.2 respectively in the case when $K = {\bf Q}_p$, taking
$\epsilon = p^{-l}$ in Theorem 1.1 and $z =   p^{-l}$ in Theorem 1.2. If each $a_i$ is sufficiently large, by scaling the results
for $f(x)$, Theorems 1.1 and 1.2 will hold for $f_a(x)$ on all of $\{x: |x_i| \leq 1$ for all $i\}$. 
If $\xi(z)$ denotes the standard character where $\xi(\sum_{k \geq k_0}b_kp^k) = e^{2\pi i ( \sum_{k = k_0}^{-1} b_k p^k)}$
 if $k_0 < 0$ and $\xi(z) = 1$
 if $k_0 \geq 0$, then since $f_a(x)$ has integer coefficients, $\xi(p^{-l}f_a(x))$ will be constant on balls of radius $p^{-l}$. Thus the integral $I(p^{-l}) = \int_{\{x:\,|x_i| \leq 1 \, \forall i\}}
\xi(p^{-l} f_a(x))\,dx$ will be the average of the $p^{ln}$ different values $\xi(p^{-l} f_a(x))$ achieves as each $x_i$
goes through the $p^l$ different balls of radius $p^{-l}$ i.e. through the different residue classes mod $p^l$.
Thus Theorem 1.2 gives Theorem 1.4.

Similarly, in Theorem 1.1 the $x$ for which $|f_a(x)| < p^{-l}$ are exactly the $x$ such that $p^l$ divides $f_a(x)$ (viewing an element of ${\bf Q}_p$ as an 
infinite series of powers of $p$). Since $f(z)$ has integer coefficients, whenever $|f_a(x)| < p^{-l}$ one will 
also have that $|f_a(x')| < p^{-l}$  for any $x'$ such that each $x_i' = x_i$ mod $p^l$. Thus whether or not $x$ 
satisfies $|f_a(x)| <p^{-l}$ depends on what each $x_i$ mod $p^l$ is. Hence the bounds on the measure of the $x$ for which 
$|f_a(x)| < p^{-l}$, given by Theorem 1.1, translate into Theorem 1.3 and we are done. 

\noindent {\bf Acknowledgement.} The author would like to thank the referee for several helpful suggestions.

\noindent {\bf 4. References.}

\noindent [AGV] V. Arnold, S Gusein-Zade, A Varchenko, {\it Singularities of differentiable maps
Volume II}, Birkhauser, Basel, 1988.\parskip = 3pt\baselineskip = 3pt

\noindent [Cl] R. Cluckers, {\it Analytic Van der Corput Lemma for $p$-adic oscillatory integrals, singular Fourier 
transforms, and restriction theorems}, Expositiones Mathematicae {\bf 29} (2011)  No. 4, 371-386.

\noindent [D] J. Denef, {\it Report on Igusa's local zeta function}, Séminaire Bourbaki, Vol. 1990/91. Astérisque No. 201-203 (1991), Exp. No. 741, 359-386 (1992). 

\noindent [DHo] J. Denef, K. Hoornaert, {\it  Newton polyhedra and Igusa's local zeta
function}, J. Number Theory {\bf 89}  (2001), no. 1, 31-64.

\noindent [DLo] J. Denef, F. Loeser, {\it Motivic Igusa zeta functions}, J. Algebraic Geom. {\bf 7} (1998), no. 3, 505-537. 

\noindent [Go] F. Gouvea, {\it $p$-adic numbers. An introduction}, Second edition. Universitext. Springer-Verlag, Berlin, 1997. vi+298 pp.

\noindent [G1] M. Greenblatt, {\it A Coordinate-dependent local resolution of singularities and 
applications},  J. Funct. Anal.  {\bf 255}  (2008), no. 8, 1957-1994.

\noindent [G2] M. Greenblatt, {\it Resolution of singularities, asymptotic expansions of oscillatory 
integrals, and related Phenomena}, J. Analyse Math. {\bf 111} (2011) no. 1, 221-245.  

\noindent [G3] M. Greenblatt, {\it Oscillatory integral decay, sublevel set growth, and the Newton
polyhedron}, Math. Annalen {\bf 346} (2010), no. 4, 857-895.

\noindent [G4] M. Greenblatt, {\it A constructive elementary method for local resolution of singularities}, submitted.

\noindent [H1] H. Hironaka, {\it Resolution of singularities of an algebraic variety over a field of characteristic zero I}, 
 Ann. of Math. (2) {\bf 79} (1964), 109-203.

\noindent [H2] H. Hironaka, {\it Resolution of singularities of an algebraic variety over a field of characteristic zero II},  
Ann. of Math. (2) {\bf 79} (1964), 205-326. 

\noindent [I1] J. Igusa, {\it Complex powers and asymptotic expansions I}, J. Reine Agnew. Math, {\bf 268/269} (1974) 110-130.

\noindent [I2] J. Igusa, {\it Complex powers and asymptotic expansions II}, J. Reine Agnew. Math, {\bf 278/279} (1975) 307-321.

\noindent [I3] J. Igusa, {\it Lectures on forms of higher degree}, Tata Inst. Fund. Research, Bombay, 1978.

\noindent [LM] B. Lichtin, D. Meuser, {\it Poles of local zeta functions and Newton polygons}, Compositio Math. {\bf 55}
 (1985), 313-332.

\noindent [R] K. Rogers, {\it A van der Corput lemma for the $p$-adic numbers}, Proc. Amer. Math. Soc.
{\bf 133} (2005), no. 12, 3525-3534.

\noindent [V] A. N. Varchenko, {\it Newton polyhedra and estimates of oscillatory integrals}, Functional 
Anal. Appl. {\bf 18} (1976), no. 3, 175-196.

\noindent [Ve] W. Veys, {\it Determination of the poles of the topological zeta function for curves}, Manuscripta Math. 
{\bf 87} (1995), no. 4, 435-448.

\noindent [Wr] J. Wright, {\it Exponential sums and polynomial congruences in two variables: the quasi-homogeneous case}, preprint.

\noindent [Zu1] W. Zuniga-Galindo, {\it Igusa's local zeta functions of semiquasihomogeneous polynomials}, Trans. Amer. Math. Soc. 353 (2001), no. 8, 3193-3207.

\noindent [Zu2] W. Zuniga-Galindo, {\it Local zeta functions and Newton polyhedra}, Nagoya Math. J. {\bf 172} (2003), 31-58. 

\line{}
\line{}

\noindent Department of Mathematics, Statistics, and Computer Science \hfill \break
\noindent University of Illinois at Chicago \hfill \break
\noindent 322 Science and Engineering Offices \hfill \break
\noindent 851 S. Morgan Street \hfill \break
\noindent Chicago, IL 60607-7045 \hfill \break
\noindent greenbla@uic.edu
\end